\documentclass[review, 3p, onecolumn]{elsarticle}

\usepackage[utf8x]{inputenc}
\usepackage{graphicx}
\usepackage{amssymb}
\usepackage{amsmath}
\usepackage{lineno}
\usepackage{subfig}
\usepackage{booktabs}
\usepackage{gensymb}
\usepackage{listings}
\usepackage{multirow}
\usepackage{lscape}
\PassOptionsToPackage{hyphens}{url}\usepackage[colorlinks=True]{hyperref}
\usepackage{algorithmic}

\biboptions{sort&compress}

\usepackage{array}
\newcolumntype{L}[1]{>{\raggedright\let\newline\\\arraybackslash\hspace{0pt}}m{#1}}
\newcolumntype{C}[1]{>{\centering\let\newline\\\arraybackslash\hspace{0pt}}m{#1}}
\newcolumntype{R}[1]{>{\raggedleft\let\newline\\\arraybackslash\hspace{0pt}}m{#1}}

\def\equationautorefname~#1\null{(#1)\null}

\usepackage{framed} 
\usepackage{multicol} 
\usepackage{nomencl} 
\makenomenclature
\renewcommand*\nompreamble{\begin{multicols}{2}}
\renewcommand*\nompostamble{\end{multicols}}
\usepackage{etoolbox}
\renewcommand\nomgroup[1]{%
  \item[\bfseries
  \ifstrequal{#1}{L}{Latin letters}{%
  \ifstrequal{#1}{M}{Greek letters}{%
  \ifstrequal{#1}{S}{Subscripts and superscripts}{}}}%
]}

\journal{Elsevier}

\begin{document}

\begin{frontmatter}

\title{Spatial analysis of thermal groundwater use based on optimal sizing and placement of well doublets}

\author[ENS]{Smajil Halilovic\corref{correspondingauthor}}
\ead{smajil.halilovic@tum.de}
\cortext[correspondingauthor]{Corresponding author}
\author[HYD,RKU]{Fabian Böttcher}
\author[HYD]{Kai Zosseder}
\author[ENS]{Thomas Hamacher}
\address[ENS]{Technical University of Munich, Chair of Renewable and Sustainable Energy Systems, Germany}
\address[HYD]{Technical University of Munich, Chair of Hydrogeology, Germany}
\address[RKU]{Department for Climate and Environmental Protection (RKU), City of Munich, Germany}

\begin{abstract}
This paper proposes an approach to optimize the technical potential of thermal groundwater use by determining the optimal sizing and placement of extraction-injection well doublets. The approach quantifies the maximum technically achievable volume of extracted groundwater in a given area and, hence, the amount of heat exchanged with the aquifer, considering relevant regulatory and hydraulic constraints. The hydraulic constraints ensure acceptable drawdown and rise of groundwater in extraction and injection wells for sustainable use, respectively, prevention of internal hydraulic breakthroughs, and adequate spacing between neighboring doublets. Analytical expressions representing these constraints are integrated into a mixed-integer linear optimization framework allowing efficient application to relatively large areas. The applicability of the approach is demonstrated by a real case study in Munich, where the geothermal potential of each city block is optimized independently. Six optimization scenarios, differing in terms of required minimum installed doublet capacity and spacings between doublets, underline the adaptability of the approach. The approach provides a comprehensive and optimized potential assessment and can be readily applied to other geographic locations. This makes it a valuable tool for thermal groundwater management and spatial energy planning, such as the planning of fourth and fifth generation district heating systems.
\end{abstract}


\end{frontmatter}

\begin{table*}[th!]
\footnotesize
\label{tab:nomenclature}
  \begin{framed}
    \nomenclature[L, B]{$B$}{Saturated aquifer thickness [m]}
    \nomenclature[L, C]{$C$}{Condition for relative well placement} 
    \nomenclature[L, N]{$N$}{Number of potential doublets}
    \nomenclature[L, S]{$S$}{Set of potential doublets}
    \nomenclature[L, Ik]{$I$}{Set of potential injection wells}
    \nomenclature[L, Ek]{$E$}{Set of potential extraction wells}
    \nomenclature[L, din]{$d_{\mathrm{inj}}$}{Decision variables for injection wells}
    \nomenclature[L, dex]{$d_{\mathrm{ext}}$}{Decision variables for extraction wells}
    \nomenclature[L, d]{$d$}{Decision variables for doublets}
    \nomenclature[L, q]{$q$}{Pump rate of a doublet [m$^3$/s]}
    \nomenclature[L, qd]{$q_{\mathrm{d}}$}{Pump rate at the drawdown threshold [m$^3$/s]} 
    \nomenclature[L, qf]{$q_{\mathrm{f}}$}{Injection rate at the upconing threshold [m$^3$/s]} 
    \nomenclature[L, qb]{$q_{\mathrm{b}}$}{Pump rate at the hydraulic breakthrough threshold [m$^3$/s]}  
    \nomenclature[L, qbl]{$q_{\mathrm{block}}$}{Pump rate per city block [m$^3$/s]}  
    \nomenclature[L, qt]{$q_{\mathrm{total}}$}{Total pump rate [m$^3$/s]}
    \nomenclature[L, qmax]{$q_{\mathrm{max}}$}{Pre-computed maximum pump rate of a doublet [m$^3$/s]}
    \nomenclature[L, qmin]{$q_{\mathrm{min}}$}{Predefined minimum pump rate of a doublet [m$^3$/s]}    
    \nomenclature[L, K]{$K$}{Hydraulic conductivity [m/s]}
    \nomenclature[L, rD]{$r_\Delta$}{External-internal well distance ratio [-]} \nomenclature[L, hn]{$h_\mathrm{n}$}{Natural groundwater level [m]} 
    \nomenclature[L, hm]{$h_{\mathrm{max}}$}{Maximum allowed groundwater level [m]}  
    \nomenclature[L, hg]{$\nabla h$}{Hydraulic gradient [-]} 
    \nomenclature[L, vD]{$v_D$}{Darcy velocity [m/s]}   
    \nomenclature[L, n]{$n_{\mathrm{doublet}}$}{Number of installed doublets in city block [-]} 
    \nomenclature[M]{$\alpha$}{Hydraulic breakthrough parameter [m$^2$/s]} 
    \nomenclature[M]{$\chi$}{Line length [m]}
    \nomenclature[M]{$\Delta$}{Distance between two wells or two doublets [m]}
    \nomenclature[M]{$\Delta_\mathrm{min}$}{Regulatory minimum distance [m]}
    \nomenclature[L, m]{$m$}{Interference parameter [m]}    
    \nomenclature[L, x]{$x, y$}{Coordinates of a well [m]}
    \nomenclature[S]{$j$}{Counter for extraction wells}
    \nomenclature[S]{$i$}{Counter for injection wells}
    \nomenclature[S]{$k, p$}{Counters for doublets}
    \nomenclature[S]{$\mathrm{ext}$}{Extraction}
    \nomenclature[S]{$\mathrm{inj}$}{Injection}
    \nomenclature[S]{$\mathrm{min}$}{Minimum}
    \nomenclature[S]{$\mathrm{max}$}{Maximum}
    \nomenclature[S]{$\tilde{a}$}{Median of $a$}
    \nomenclature[L, u]{$\boldsymbol{u}$}{Optimization variables}
    \nomenclature[L, g1]{$\mathbf{g}_1$}{Equality constraints}
    \nomenclature[L, g2]{$\mathbf{g}_2$}{Inequality constraints}    
    \printnomenclature
  \end{framed}
\end{table*}

\section{Introduction}\label{sec:Introduction}

Shallow geothermal energy (SGE) plays an increasingly important role in the decarbonization of the heating and cooling sector \citep{Lund.2021}. 
Especially in the context of 4th generation district heating (4GDH) systems \citep{Lund.2014} and 5th generation district heating and cooling (5GDHC) systems \citep{Boesten.2019}, SGE systems represent a promising technology as their application is expanded from individual users to communities in this case. 
One way of SGE utilization is through open-loop systems, commonly referred to as groundwater heat pumps (GWHPs). 
These systems directly exploit the thermal energy of groundwater through extraction-injection well doublets by pumping groundwater from extraction wells and returning it to the same shallow aquifer through injection wells after thermal exchange \citep{Stauffer.2014}. 
Consequently, the properties of groundwater, including its quantity, quality, depth, and temperature, are the most important factors affecting the feasibility and performance of GWHP systems \citep{Banks.2012, Lee.2006}.
Accurate characterization and consideration of these groundwater properties is essential for ensuring sustainable and efficient operation of GWHPs \citep{Busby.2009}.

In the context of groundwater utilization for GWHP systems, it is crucial to recognize that groundwater availability and properties exhibit significant spatial variation \citep{Fry.2009}.
In addition, spatial availability for GWHP well installations in urban areas is limited due to extensive building development.
Therefore, conducting a spatial analysis is essential to identify adequate well locations and sizing for well doublets. 
The goal of such an analysis is to quantify suitable groundwater extraction values for the targeted urban energy planning level, such as a plot of land or city block \citep{Schiel.2016}. 
Accurate estimation of groundwater extraction and injection rates, can ensure sustainable GWHP operation based on the local hydrogeological conditions and in compliance with relevant legal and technical constraints. 
The technical potential derived from such an analysis is a basis for active thermal groundwater management and the development of urban energy strategies \citep{GarciaGil.2022}, including the planning of 4GDH and 5GDHC systems.

Several research studies assess the potential of thermal groundwater use at different locations and considering various constraints \citep{Epting.2018, Epting.2020, Casasso.2017, Munoz.2015, GarciaGil.2015, Pujol.2015, Arola.2014, Arola.2014b, Bottcher.2019}. 
Some studies focus on specific technical and/or regulatory constraints, while others combine multiple constraints to estimate the technical potential, i.e. technically feasible groundwater pumping rates, at a given spatial resolution. 
In general, their aim is to estimate the local geothermal potential, but not to optimize the technical potential.
However, the geothermal potential can be maximized through strategic sizing and placement of well doublets within the considered area. 
As the thermal use of groundwater with well doublets induces hydraulic and thermal changes in the aquifer, each operating doublet consumes space and obstructs the installation of additional doublets. 
In the vicinity of wells, pumping may create a considerable drawdown and injection an upconing of groundwater, respectively. 
Especially in the planning stage of larger GWHP systems with multiple well doublets, a consideration of hydraulic influences is crucial for a sustainable well arrangement, as each well doublet leaves its own hydraulic footprint in the aquifer.
These characteristics of multi-doublet systems have not been addressed in the existing potential assessment studies.
This aspect is particularly important when the spatial planning level allows flexibility in the arrangement of doublets and associated wells. 
For instance, analyzing the geothermal potential for large GWHPs in the context of future 4GDH and 5GDHC systems requires determining the optimal size and placement of multiple wells simultaneously.
To answer the question of optimal well count, sizing and placement in the potential assessment, it is necessary to integrate optimization methods into the analysis.

\citet{Halilovic.2023a} recently reviewed optimization approaches for GWHP systems and concluded that there are only a few studies addressing the topic of optimal well placement and/or sizing. 
Some of the existing optimization approaches employ numerical groundwater simulation models and are therefore not suitable for integration into large-scale potential assessments due to their high computational costs \citep{Park.2021, Park.2020, Halilovic.2022a}. 
In contrast, optimization approaches based on analytical models prove to be more suitable for such integration, as they generally require less computational cost while providing reasonably accurate and conservative estimates. 
To date, only one research study has implemented an optimization approach based on an analytical groundwater simulation model \citep{Halilovic.2023b}. 
The approach uses the linear advective heat transport model (LAHM) to estimate thermal plumes caused by GWHPs \citep{Kinzelbach.1987}. 
The developed approach optimally places GWHPs and their associated wells in the considered area in order to minimize thermal interactions between the systems and simultaneously maximize the heat extracted from the groundwater, thereby maximizing the spatial potential of thermal groundwater use. 
This optimization approach is promising for potential estimation studies and is already applied on city-scale supporting the municipal heat planning of Munich \citep{Halilovic.2023b}, but has certain limitations.

The main limitations of the approach proposed in \citet{Halilovic.2023b} arise from the characteristics of the LAHM model, which assumes homogeneous groundwater conditions throughout the entire study area. 
Additionally, the approach does not consider any hydraulic aspects, such as pumping limits due to induced groundwater drawdowns in extraction wells and the resulting hydraulic footprint that prevents additional wells of the same type to be installed nearby.
Furthermore, systems' pumping rates are predefined based on the estimated energy demand of the corresponding plots, which means that only the placement of the systems and their wells is optimized and their sizing remains unchanged. 
This makes the approach unsuitable for certain applications, such as potential analysis for large multi-doublet GWHPs in 4GDH and 5GDHC systems, since hydraulic constraints are crucial in this case. 
Therefore, there is a need for an optimization approach that accounts for aspects not covered in \citet{Halilovic.2023b}, particularly hydraulic constraints and spatial heterogeneity of groundwater parameters.

The main objective of this paper is to introduce a novel approach for optimizing the technical potential of thermal groundwater use. 
The proposed approach simultaneously optimizes the sizing and placement of doublets and associated wells within feasible areas, with the goal of maximizing the geothermal potential, i.e. the extracted heat from groundwater. 
To achieve this, the approach considers multiple important factors, including regulatory constraints, spatial heterogeneity of groundwater properties, and relevant hydraulic constraints. 
The latter includes considerations of drawdown in extraction wells, groundwater rise in injection wells, internal hydraulic breakthroughs, and spacing between adjacent well doublets based on their hydraulic footprints \citep{Clyde.1983}. 
To effectively include these constraints, the method integrates analytical expressions for pumping rate limits from \citet{Bottcher.2019} into a mixed-integer linear optimization framework. 
The new approach provides a comprehensive and robust potential estimation, as demonstrated by using a real case study wherein the technical potential of thermal groundwater use of each city block of Munich is optimized.

The paper is organized as follows: Section~\ref{sec:methodology} describes the methodology, followed by its implementation and a case study in Section~\ref{sec:Implementation}. Section~\ref{sec:Results} presents and analyzes the results. Section~\ref{sec:discussion} discusses the advantages, limitations and possible future improvements and applications of the approach. The paper concludes with a summary in Section~\ref{sec:Summary}.

\section{Methodology}
\label{sec:methodology}

The proposed approach combines analytical expressions that describe groundwater pumping limits (Section~\ref{ssec:TAP}) with mixed-integer linear programming techniques (Section~\ref{ssec:optimization}) to optimize the placement and sizing of extraction-injection well doublets.

\subsection{The TAP method}
\label{ssec:TAP}

\citet{Bottcher.2019} developed the Thermal Aquifer Potential (TAP) method to analyze the technical potential of thermal groundwater use.
This method is based on empirical analytic formulas that describe the maximum pumping rates of a well doublet with respect to three different constraints: maximum drawdown in the extraction well, upconing threshold at the injection well, and hydraulic breakthrough between the two wells.

The TAP method estimates the maximum pumping rate of a doublet at the drawdown threshold $q_{\mathrm{d}}$ as follows:
\begin{equation}
q_{\mathrm{d}} = 0.195\cdot K\cdot B^2\, ,  \label{eq:tap_drawdown}
\end{equation}
\noindent where $K$ is the hydraulic conductivity and $B$ is the saturated aquifer thickness. 
The considered threshold for drawdown is one-third of the saturated aquifer thickness \citep{Bottcher.2019}.
The maximum injection rate at the upconing (flooding) threshold $q_{\mathrm{f}}$ is estimated with:
\begin{equation}
q_{\mathrm{f}} = (h_{\mathrm{max}}-h_\mathrm{n})\cdot K\cdot B^{0.798}\cdot\exp(29.9\cdot \nabla h)\, ,  \label{eq:tap_flood}
\end{equation}
\noindent where $h_{\mathrm{max}}$ and $h_\mathrm{n}$ are the maximum allowed and the natural groundwater level, respectively, and $\nabla h$ is the hydraulic gradient.
Finally, the TAP method calculates the maximum pumping rate at the hydraulic breakthrough threshold $q_{\mathrm{b}}$ of the well doublet using the following equation:
\begin{equation}
q_{\mathrm{b}} = \frac{\pi}{1.96}\cdot v_D\cdot B\cdot \, \Delta_\mathrm{wells},  \label{eq:tap_break}
\end{equation}
\noindent where $v_D$ is the Darcy velocity and $\Delta_\mathrm{wells}$ is the internal distance between the extraction and injection well of the doublet.
In this work, to simplify the description of the optimization problem in later sections, \autoref{eq:tap_break} is rewritten as follows:
\begin{equation}
q_{\mathrm{b}} = \alpha \cdot \, \Delta_\mathrm{wells},  \label{eq:tap_break_2}
\end{equation}
\noindent where $\alpha = (\pi/1.96)\cdot v_D\cdot B$ is the hydraulic breakthrough parameter of the considered well doublet.
The TAP method defines the technical pumping rate of a well doublet as the minimum of the three previously specified pumping rates $q_{\mathrm{d}}$, $q_{\mathrm{f}}$ and $q_{\mathrm{b}}$.

\begin{figure}[tb]
\centering
\includegraphics[width=0.4\textwidth]{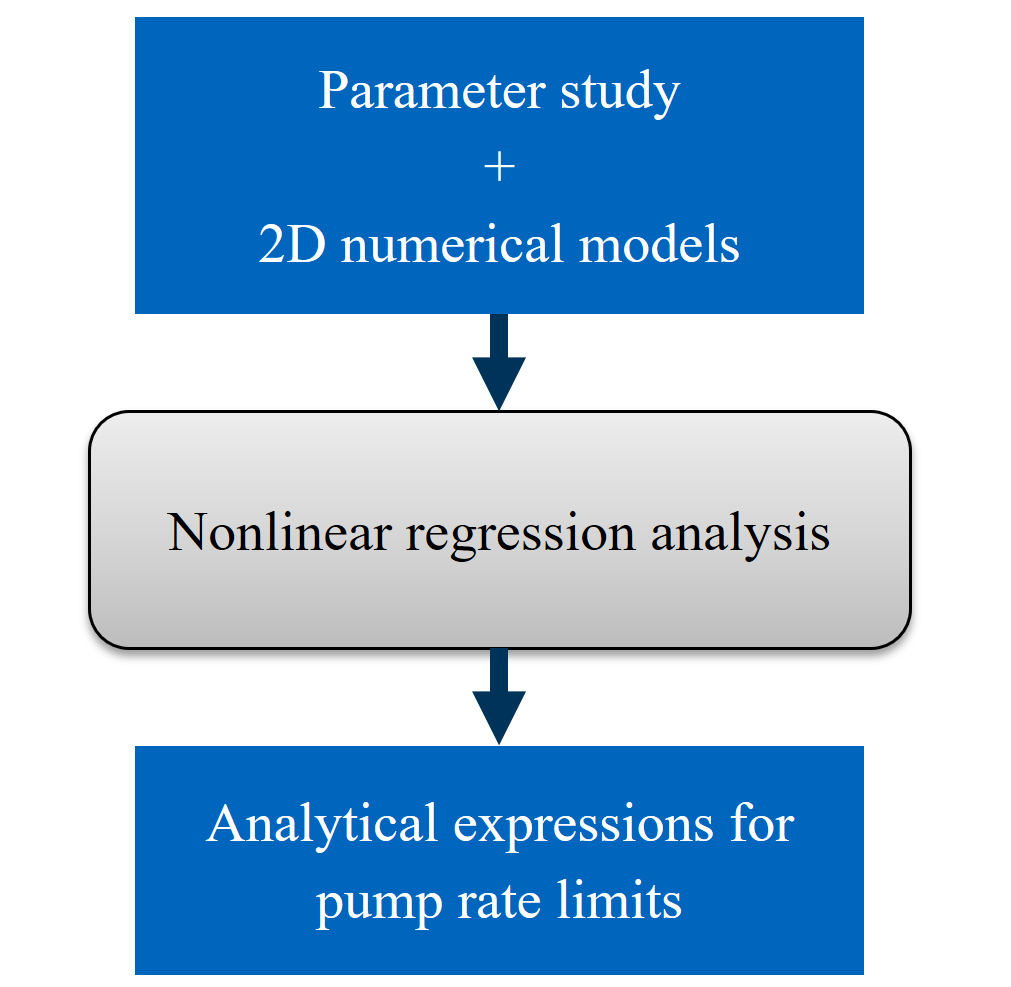}
\caption{Flowchart of the TAP method.}
\label{fig:tap_method}
\end{figure}

The analytical formulas \autoref{eq:tap_drawdown}-\autoref{eq:tap_break} in the TAP method are derived from the results of numerical parameter studies using nonlinear regression analysis.
\autoref{fig:tap_method} shows the overall flowchart of the TAP method.
In the first phase, idealized 2D box models are prepared for numerical groundwater simulations in the parameter study. 
In this study, important parameters are varied within a reasonable range of values, and steady-state simulations are performed to gain conservative results. 
The extraction and injection wells are placed in the center of the 2D models parallel to the groundwater flow direction and at different distances from each other. 
Subsequently, the results of the numerical simulations from the parameter study are used in a nonlinear regression analysis to fit the analytical expressions \autoref{eq:tap_drawdown}-\autoref{eq:tap_break}.

In addition to the threshold pumping rates $q_{\mathrm{d}}$, $q_{\mathrm{f}}$, and $q_{\mathrm{b}}$, the authors of the TAP method also analyzed the hydraulic footprint of a well doublet using idealized 2D models in a similar manner.
This analysis is necessary for estimating technical potential because neighboring systems limit the available water budget and therefore increase the likelihood of hydraulic breakthrough within the system if they are located too close.
Thus, the authors determined a correlation between the percentage of cycled water (inter-flow) in a well doublet and the external-internal well distance ratio $r_\Delta$, i.e. the ratio between the distance to neighboring doublets and the inner well spacing.

\subsection{Optimization}
\label{ssec:optimization}

An optimization procedure can be used to maximize the technical potential of thermal groundwater use through optimal placement and sizing of well doublets.
The proposed optimization approach integrates the equations from the TAP method into a mixed-integer linear program to maximize the thermal aquifer potential while satisfying technical and legal constraints.
The placement of doublets and their wells is based on a selection of predefined potential well locations, which are defined in a pre-processing step as described in the following.

\subsubsection{Definition of potential well locations and doublets}
\label{sssec:positions}

\begin{figure}[tb]
\centering
\captionsetup[subfigure]{position=top}
\subfloat[Potential well locations]{\label{fig:Potential_wells}\includegraphics[width=0.52\textwidth]{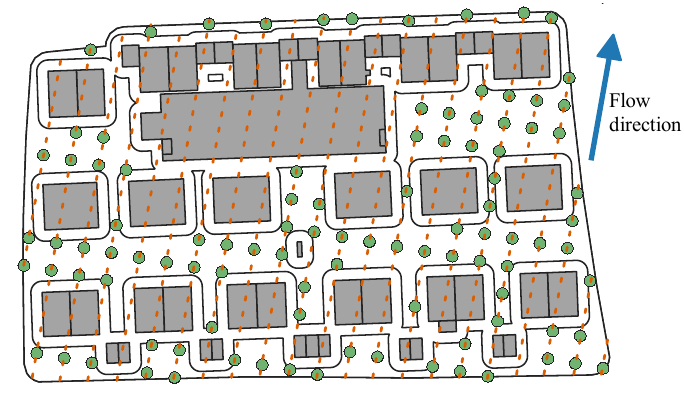}}
\subfloat[Optimal well locations]{\label{fig:Selected_wells}\includegraphics[width=0.48\textwidth]{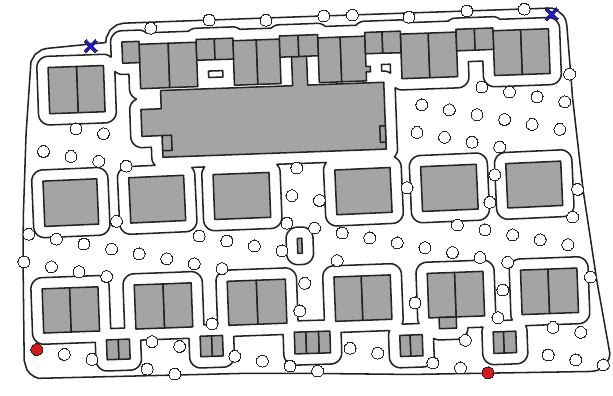}}
\\
\vspace{-.5cm}
\hspace{1cm}\subfloat{\includegraphics[width=0.24\textwidth]{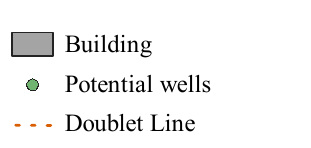}}\hspace{4cm}
\subfloat{\includegraphics[width=0.25\textwidth]{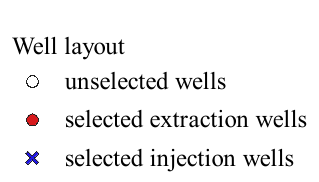}}
\vspace{-.3cm}
\caption{Selection of optimal well locations: (a) Potential wells. (b) Optimal (selected) wells and corresponding doublets.
\label{fig:wells_selection}}
\end{figure}

In the pre-processing phase, potential well locations are determined taking into account ground plans of the existing buildings, legal constraints and groundwater flow direction in the considered area. 
To determine a feasible area for well placement, a minimum distance of 3 meters between wells and buildings is maintained in the first step, by applying corresponding buffer zones (see Figure~\ref{fig:Potential_wells}). 
Within this delineated area, potential wells are strategically positioned at the nodes of a virtual grid that is oriented according to the groundwater flow direction at the centroid of the corresponding polygon. 
The grid is designed with a constant well-to-well spacing, with one of its axes aligned with the groundwater flow direction and the other axis perpendicular to it. 
The grid is then divided into lines parallel to groundwater flow, which are denoted as doublet lines in Figure~\ref{fig:Potential_wells}, and wells are grouped based on the lines they lie on. 
The simplifying assumption of this procedure is that only wells placed on the same line are allowed to be installed as an extraction-injection well doublet for the potential multi-well system. 
This ensures that the installed wells are aligned with the groundwater flow direction, which is a prerequisite for applying the TAP method (see Section~\ref{ssec:TAP}). The pre-processing step results in hydrogeologically and legally viable potential well locations, which are further used in the optimization procedure.

\subsubsection{Optimization concept}
\label{sssec:well_selection}

As described previously, each line in Figure~\ref{fig:Potential_wells} corresponds to one potential well doublet, i.e. one upstream extraction and one downstream injection well which are installed according to the groundwater flow direction.
There are multiple potential well locations on each line and multiple potential doublets (lines) in a city block, which gives a high degree of freedom in the layout (installation) of a large multi-doublet system and the placement of its wells.
In addition, the size of a doublet, i.e. its pumping rate, is interdependent with the well locations due to the hydraulic breakthrough limits. 
Therefore, an optimization procedure is required to determine the optimal combination of doublets to be installed in a city block, i.e. the size and placement of the doublets and the placement of the corresponding wells.
The optimization concept is described below:
\begin{itemize}
    \item the area of interest has a predefined number of potential well doublets (lines),
    \item each doublet (line) has a predefined number of potential well locations, which can be extraction or injection wells,
    \item for each potential doublet (line) there are two optimization variables: the binary variable $d$ corresponding to the decision whether to install ($d=1$) or not ($d=0$) this doublet, and the continuous variable $q\in\mathbb{R}^{+}_0$ representing the pumping rate of the doublet, 
    \item for each potential well location, there are two binary optimization variables: $d_{\mathrm{ext}}$ and $d_{\mathrm{inj}}$, which represent the selection decision for the extraction and injection well, respectively,
    \item if an extraction well is installed (selected) at the considered well location, then $d_{\mathrm{ext}}=1$ and $d_{\mathrm{inj}}=0$, and vice versa, in case of an injection well, $d_{\mathrm{ext}}=0$ and $d_{\mathrm{inj}}=1$,
    \item if neither an extraction nor an injection well is selected at the considered potential well location, then $d_{\mathrm{ext}}=d_{\mathrm{inj}}=0$,
    \item if the well doublet is installed ($d=1$), one extraction well and one injection well is selected from the potential well locations on the corresponding line,
    \item if doublet is not installed ($d=0$), all potential wells on that line are deselected.
\end{itemize}
Figure~\autoref{fig:Selected_wells} shows an example of a city block with two installed doublets and the placement (selection) of their wells.
It is important to point out that the proposed approach does not consider thermal interactions between doubles, specifically the propagation of thermal plumes in groundwater induced by GWHPs.
Therefore, the approach focuses on optimizing the technical potential of smaller areas (e.g. city blocks) separately and does not address the optimization of the spatial potential of multiple areas jointly.
In the next section, further details on the relations between the optimization variables described previously, as well as on other optimization constraints, are provided.

\subsubsection{Objective function and constraints}
\label{sssec:optim_problem}

The optimization objective is to maximize the technical potential of thermal aquifer utilization in a city block, i.e. to maximize groundwater extraction by well doublets while meeting the corresponding technical and legal constraints. 
Thus, the objective function to be maximized is defined as:
\begin{equation}
q_{\mathrm{total}} = \sum_{k=1}^{N} q_k\,,  \label{eq:obj1}
\end{equation}
\noindent where $q_{\mathrm{total}}$ is the total pumping rate of all doublets, $q_k$ is the pumping rate of a single doublet $k$ and $N$ is the number of potential doublets.
The maximization of the volume of extracted groundwater simultaneously maximizes the thermal energy exchange with the aquifer.

The problem also includes several optimization constraints related to the installation and operation of doublets and their wells.
The first set of constraints ensures that only installed doublets can be operated:
\begin{equation}
q_k \leq q_{\mathrm{max},k}\cdot d_k \quad \forall k\in S \,, \label{eq:cst_operation}
\end{equation}
\noindent where $q_{\mathrm{max},k}$ is the pre-computed maximum pumping rate of a doublet $k$ and $S=\{1, ..., N\}$ is the set of all potential well doublets.
If a doublet is not installed ($d_k=0$), it cannot pump groundwater ($q_k=0$), otherwise its pumping rate is limited by $q_{\mathrm{max},k}$, which is calculated in the pre-processing as follows:
\begin{equation}
q_{\mathrm{max},k} = \min(\max\limits_{j\in E_k} q_{\mathrm{d},j}, \max\limits_{i\in I_k} q_{\mathrm{f},i}) \,, \label{eq:cst_q_max}
\end{equation}
\noindent where $q_{\mathrm{d},j}$ and $q_{\mathrm{f},i}$ are the "threshold" pumping and injection rates from the TAP method, i.e. equations \autoref{eq:tap_drawdown} and \autoref{eq:tap_flood}, respectively, and $E_k$ and $I_k$ are the sets of all potential extraction and injection wells of the doublet $k$, respectively.
Thus, $q_{\mathrm{max},k}$ represents the minimum between the two: the maximum pumping rate of all potential extraction wells $j$, based on the drawdown threshold, and the maximum injection rate of all potential injection wells $i$, based on the upconing threshold, of the doublet $k$ (see Section~\ref{ssec:TAP}).
This initial estimate of the theoretical upper bound for the pumping rates is used in the optimization constraints to reduce the exploratory design space and thereby speed up the overall optimization process.

The second set of constraints specifies a minimum pumping rate for installed well doublets: 
\begin{equation}
q_\mathrm{min}\cdot d_k \leq q_k \quad \forall k\in S \,, \label{eq:cst_q_min}
\end{equation}
where $q_\mathrm{min}$ is the predefined minimum pumping rate of a doublet.
These constraints serve to prevent the installation of too small doublets, which are not economically viable in practice.

The third set of constraints corresponds to the fact that each doublet consists of a single extraction-injection well pair:
\begin{equation}
\sum_{i\in I_k} d_{\mathrm{inj},i} = d_k \quad \forall k\in S \,, \label{eq:cst1}
\end{equation}
\begin{equation}
\sum_{j\in E_k} d_{\mathrm{ext},j} = d_k \quad \forall k\in S\,. \label{eq:cst2}
\end{equation}
\noindent This also implies that the number of installed extraction and injection wells of the same doublet must be the same, i.e. 0 or 1, depending on whether the doublet is installed or not.

The next set of constraints are limitations on the pumping rates of well doublets based on the TAP method. 
The first group of such constraints ensures that none of the installed extraction wells (doublets) exceeds the drawdown threshold defined in the TAP method:
\begin{equation}
q_k \leq d_{\mathrm{ext},j}\cdot q_{\mathrm{d},j} + q_{\mathrm{max},k}\cdot(1-d_{\mathrm{ext},j}) \quad \forall j\in E_k\,, \forall k\in S \,. \label{eq:cst_drawdown}
\end{equation}
\noindent Depending on which extraction well $j$ is selected ($d_{\mathrm{ext},j}=1$), the pumping rate $q_k$ of the doublet is limited by the pumping rate at the drawdown threshold of this well $q_{\mathrm{d},j}$.
If the well $j$ is not selected ($d_{\mathrm{ext},j}=0$), the constraint reads as $q_k \leq q_{\mathrm{max},k}$, which should hold in any case.
The form of the constraint \autoref{eq:cst_drawdown} is necessary to ensure that only the selected extraction wells set an upper limit on the pumping rate based on the drawdown threshold.
The other constraints in connection with the TAP method are formulated in a similar manner.

The second group of TAP-related constraints limits the injection rates based on the upconing threshold:
\begin{equation}
q_k \leq d_{\mathrm{inj},i}\cdot q_{\mathrm{f},i} + q_{\mathrm{max},k}\cdot(1-d_{\mathrm{inj},i}) \quad \forall i\in I_k\,, \forall k\in S \,. \label{eq:cst_flood}
\end{equation}
\noindent Similar to \autoref{eq:cst_drawdown}, these constraints enforce that the selected injection wells, and thus the corresponding doublets, do not exceed the upconing threshold. 

The next group of TAP-related constraints ensures that an internal hydraulic breakthrough is prevented by limiting the doublet's pumping rate:
\begin{equation}
q_k \leq \alpha_{i,j}\cdot \Delta_{i,j} + q_{\mathrm{max},k}\cdot(2-d_{\mathrm{ext},j}-d_{\mathrm{inj},i}) \quad \forall j\in E_k\,, \forall i\in I_k\,, \forall k\in S \,, \label{eq:cst_int_break}
\end{equation}
\noindent where $\Delta_{i,j}$ is the distance between injection and extraction wells $i$ and $j$, respectively, and $\alpha_{i,j}=(\alpha_{i}+\alpha_{j})/2$ is the averaged hydraulic breakthrough parameter $\alpha$ of these two wells.
The constraint \autoref{eq:cst_int_break} is activated, i.e. it becomes $q_k \leq \alpha_{i,j}\cdot \Delta_{i,j}$, only for the selected $(i,j)$ well pair ($d_{\mathrm{ext},j}=d_{\mathrm{inj},i}=1$).
In all other cases, the constraint does not introduce additional, more stringent upper limits on the pumping rate $q_k$.

In addition to the relation between pumping rate and internal well distance defined by the constraint \autoref{eq:cst_int_break}, the well placement must also comply with the regulatory minimum internal distance and the natural order (upstream-downstream) of well placement:
\begin{equation}
 d_{\mathrm{inj},i} + d_{\mathrm{ext},j} \leq 1 \quad\text{if}\quad (\Delta_{i,j}<\Delta_\mathrm{min} \; \text{or}\; \neg C_{i,j}) \quad \forall j\in E_k\,, \forall i\in I_k\,, \forall k\in S \,, \label{eq:cst_Dint_min}
\end{equation}
\noindent where $\Delta_\mathrm{min}$ is the defined regulatory minimum distance and $C_{i,j}$ is the relative well placement condition for the $(i,j)$ well pair, which states that the injection well should be placed downstream relative to the extraction well.
If the internal distance of the considered well pair $(i,j)$ is smaller than $\Delta_\mathrm{min}$ or if the condition $C_{i,j}$ is not satisfied, these two wells cannot be installed together as a doublet.
The regulatory distance $\Delta_\mathrm{min}$ between extraction and injection wells of the same doublet is defined to avoid hydraulic and thermal breakthroughs within the system \citep{Banks.2009}.
This distance is 10 m in the case study (see \autoref{ssec:case_study}), which is located in the German state of Bavaria \citep{LfU.2012}.

The remaining optimization constraints address the spacing and sizing of neighboring doublets. 
The first group of such constraints guarantees sufficient distance between two neighboring doublets considering their hydraulic footprint derived from the TAP method. 
The pumping rates of two neighboring doublets $k$ and $p$ are limited based on their mutual distance $\Delta_{k,p}$ as follows:
\begin{equation}
\frac{q_k}{\tilde{\alpha}_k} + \frac{q_p}{\tilde{\alpha}_p} \leq \frac{2}{r_\Delta}\cdot\Delta_{k,p} + m_{k,p}\cdot(2-d_k-d_p) \quad\text{if}\quad \chi_k + \chi_p > \frac{2}{r_\Delta}\cdot\Delta_{k,p} \quad \forall k,p\in S \,, \label{eq:cst_ext_break}
\end{equation}
\noindent where: $\tilde{\alpha}_k$ and $\tilde{\alpha}_p$ are the hydraulic breakthrough parameters for the doublets $k$ and $p$, respectively; $m_{k,p} = q_{\mathrm{max},k}/\tilde{\alpha}_k + q_{\mathrm{max},p}/\tilde{\alpha}_p$ is the interference parameter between these two doublets; $\chi_k$ and $\chi_p$ are the line lengths representing the doublets $k$ and $p$, respectively, i.e. the maximum possible internal well distances for the doublets; and $r_\Delta$ is the previously defined external-internal well distance ratio.   
The hydraulic breakthrough parameter $\tilde{\alpha}$ for each doublet is defined as the median value of this parameter among all potential wells within that doublet. 
To reduce computational complexity and avoid excessive constraints, the constraint is applied only to pairs of potential doublets that are relatively close to each other, since mutual hydraulic influence is relevant in this case.
The relative closeness of two doublets is determined by comparing their distance $\Delta_{k,p}$ with the averaged maximum possible internal well spacing of the doublets $(\chi_k + \chi_p)/2$, multiplied by the chosen external-internal spacing ratio $r_\Delta$. 
If $\Delta_{k,p}\geq r_\Delta\cdot(\chi_k + \chi_p)/2$, the doublets are sufficiently far apart and the constraint \autoref{eq:cst_ext_break} is not applied. 
Otherwise, the constraint is activated.

When the doublet pair $(k,p)$ is installed, i.e. $d_k=d_p=1$, the constraint \autoref{eq:cst_ext_break} takes the following form:
\begin{equation}
\frac{q_k}{\tilde{\alpha}_k} + \frac{q_p}{\tilde{\alpha}_p} \leq \frac{2}{r_\Delta}\cdot\Delta_{k,p} \,, \label{eq:cst_ext_break2}
\end{equation}
which imposes a limit on the weighted sum of the pumping rates of the doublets based on their distance.
The formulation of \autoref{eq:cst_ext_break2} is derived from the constraint for the internal hydraulic breakthrough of a doublet \autoref{eq:cst_int_break} and using the definition of the ratio $r_\Delta$.
In all other cases, i.e. when one or both doublets are not installed, the parameter $m_{k,p}$ ensures that the constraint \autoref{eq:cst_ext_break} is always satisfied, thus avoiding the introduction of any additional limitations on the pumping rates.

The second set of constraints for neighboring doublets specifies a minimum distance between two installed doublets, which is determined by the regulatory internal distance $\Delta_\mathrm{min}$ between wells within the same doublet:
\begin{equation}
 d_k + d_p \leq 1 \quad\text{if}\quad \Delta_{k,p}<r_\Delta \cdot\Delta_\mathrm{min} \quad \forall k,p\in S \,. \label{eq:cst_Dext_min}
\end{equation}
If the distance between two potential doublets is smaller than $r_\Delta \cdot\Delta_\mathrm{min}$, these doublets cannot be installed jointly. 
This constraint guarantees that neighboring doublets maintain an adequate spacing even for smaller doublet sizes (pumping rates), which is not addressed by the constraint \autoref{eq:cst_ext_break}.

Based on the previously defined objective function and constraints, the underlying optimization problem can be formulated as follows, 
\begin{subequations}
\begin{alignat}{2}
&\!\max_{\boldsymbol{u}} &\qquad& q_{\mathrm{total}}(\boldsymbol{u})\label{eq:optProb}\\
&\text{subject to} &      & \mathbf{g}_1(\boldsymbol{u})=0,\\
& &\qquad& \mathbf{g}_2(\boldsymbol{u})\leq 0,
\end{alignat}\label{eq:Prob}
\end{subequations}
\noindent where: $\boldsymbol{u}$ represents the vector of all optimization variables, which includes the binary variables $\boldsymbol{d}_{\mathrm{ext}}$, $\boldsymbol{d}_{\mathrm{inj}}$, $\boldsymbol{d}$ and the continuous variables $\boldsymbol{q}$; $\mathbf{g}_1$ represents all the equality constraints defined by \autoref{eq:cst1} and \autoref{eq:cst2}; $\mathbf{g}_2$ represents all the inequality constraints defined by \autoref{eq:cst_operation}, \autoref{eq:cst_q_min}, \autoref{eq:cst_drawdown}, \autoref{eq:cst_flood}, \autoref{eq:cst_int_break}, \autoref{eq:cst_Dint_min}, \autoref{eq:cst_ext_break} and \autoref{eq:cst_Dext_min}.
The optimization problem \autoref{eq:Prob} is a mixed-integer linear program that is solved independently for each city block, as described in the following.

\section{Implementation}\label{sec:Implementation}

The introduced optimization approach is implemented using Python-MIP \citep{Santos.2020}, an open-source package specifically designed for modeling and solving mixed-integer linear programs.
The Python code, including a functional example, is freely available at \citep{Halilovic.2023c}.
Pre-processing of potential well locations and pumping rate limits is conducted in Python using geopandas and dependent libraries \citep{kelsey_jordahl_2020_3946761}.
In the following sections, the case study area and the considered optimization scenarios are described.

\subsection{Case study}
\label{ssec:case_study}

\begin{figure}[tb]
\centering
\includegraphics[width=0.78\textwidth]{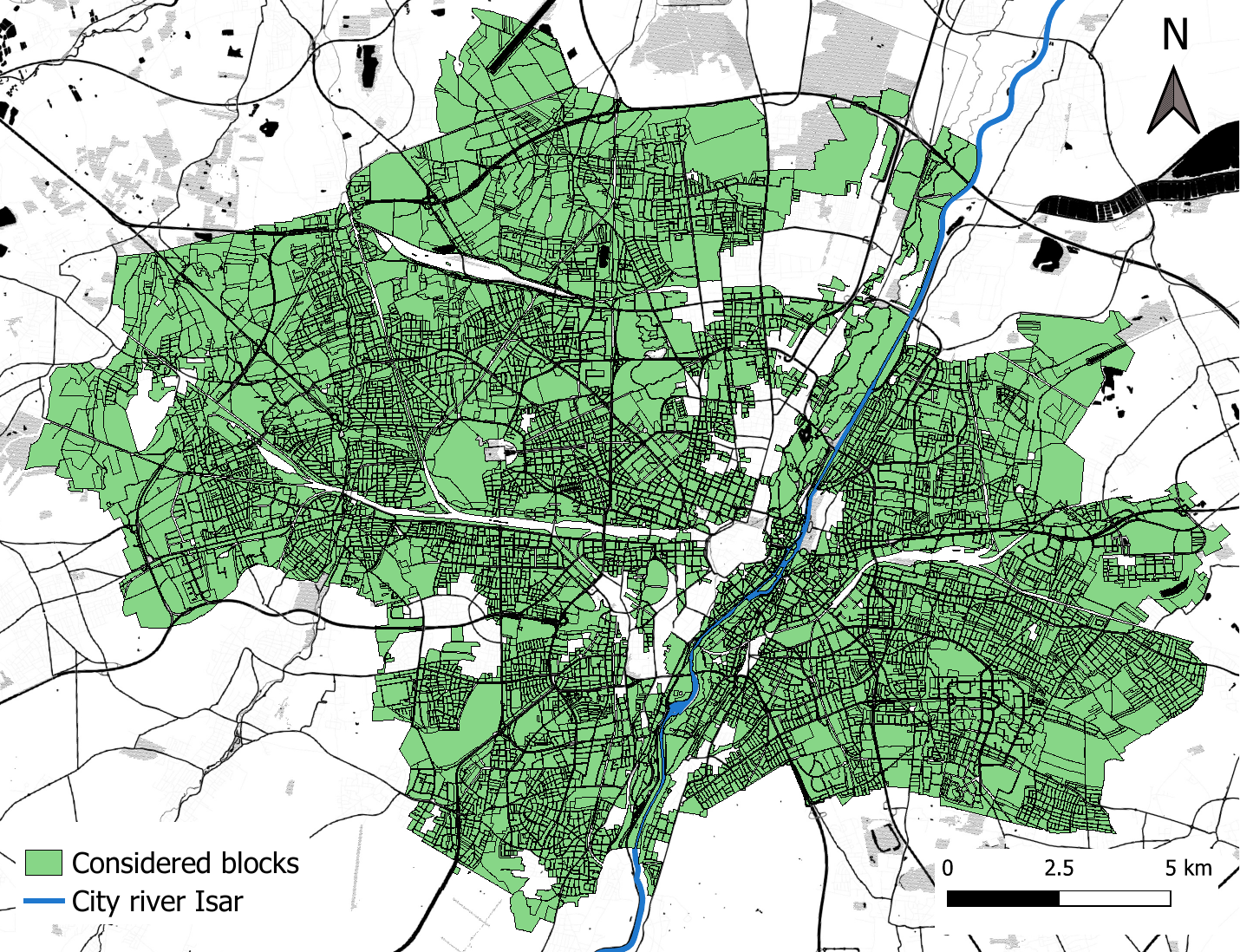}
\caption{Case study with considered city blocks.}
\label{fig:case_study}
\end{figure}

The presented optimization approach is applied on a city block level in the entire city area of Munich (see \autoref{fig:case_study}). 
This means that the optimization problem \autoref{eq:Prob} is formulated and solved for each city block, thereby optimizing the potential of each block independently from other blocks.
The city block level is selected because it is aligned with the energy planning scale relevant to future 4GDH and 5GDHC systems.
Munich offers favorable conditions for exploiting thermal energy from groundwater due to its location on a productive and shallow gravel aquifer. 
Extensive studies have been conducted to characterize key hydro-geological parameters, such as hydraulic conductivity, aquifer thickness, and groundwater flow direction, within the city area. 
Detailed information on these parameters are provided by \citet{Bottcher.2019} and \citet{Zosseder.2022}.

In the pre-processing step (see Section~\ref{sssec:positions}), potential well locations are determined for each city block by initially using a constant distance of 5 m between wells. 
Due to the large area of some city blocks, this results in a high number of potential well locations in those blocks. 
Thus, to simplify the calculation and reduce computational time, an iterative approach is used to reduce the number of potential wells per city block. 
The constant distance between wells is iteratively increased by 2.5 meters for each block until the number of potential wells within the block is reduced to 100 or less.
Once the potential well locations are determined, the relevant values from the TAP method $q_{\mathrm{d}}$, $q_{\mathrm{f}}$, and $\alpha$ are calculated for each location using the available groundwater parameter data for the city of Munich. 
In addition, a filtering process is used during pre-processing to exclude potential well locations with limited potential. 
Specifically, potential wells with pre-calculated pumping rates at drawdown or upconing thresholds ($q_{\mathrm{d}}$ or $q_{\mathrm{f}}$) below 1 [l/s] are removed. 
The value of 1 [l/s] is used in this work because the focus here is on large multi-well systems that can be used in 4GDH or 5GDHC systems.
This process ensures that areas (blocks) lacking sufficient potential due to groundwater conditions are excluded from further analysis. 
The potential analysis (optimization) is then performed for the remaining 8751 city blocks of Munich, as shown in \autoref{fig:case_study}.

\subsection{Optimization scenarios}
\label{ssec:Scenarios}

In total, six distinct optimization scenarios were investigated, each representing different combinations of two parameters: the external-internal well distance ratio $r_\Delta$ and the predefined minimum pumping rate of an installed doublet $q_{\mathrm{min}}$.
Three values were considered for $r_\Delta$: 1.5, 2, and 3. 
These values correspond to about 10\%, 5\%, and 2.5\%, respectively, of the inter-flow according to the TAP method (see Section~\ref{ssec:TAP}). 
The last case, with $r_\Delta=3$, is the most conservative and is characterized by the lowest level of interaction between neighboring doublets.
These $r_\Delta$ values were paired with two values for $q_{\mathrm{min}}$: 1 and 5 [l/s]. 
In the latter case, the use of larger doublets, e.g. for 4GDH and 5GDHC grids, is particularly emphasized.

\section{Results}\label{sec:Results}

The results obtained from the optimization scenarios are presented and analyzed in this section.
Table~\ref{tab:results_01} provides a summary of the results for all optimization scenarios introduced in Section~\ref{ssec:Scenarios}.
The table includes the following results: the total number of installed well doublets in the city of Munich, the maximum and mean number of installed doublets per city block, the average of pumping rates of the largest installed doublets in city blocks (i.e. the average of maximum pumping rates of installed doublets per city block), the mean value of pumping rates of all installed doublets (excluding non-installed doublets), the number of city blocks with and without installed doublets, the total installed pumping rate for the entire city, and the maximum and mean values of installed pumping rates per city block.

It is evident that scenarios with larger values of $r_\Delta$ and the same $q_{\mathrm{min}}$ have a smaller number of installed well doublets with larger capacities (pumping rates).
At the same time, as the $r_\Delta$ value increases, the installed capacities per city block decrease.
This result conforms with expectations, as these scenarios follow a more conservative approach that requires larger distances between neighboring doublets.
Furthermore, for a given $r_\Delta$, the maximum installed pumping rate per city block remains unchanged regardless of the $q_{\mathrm{min}}$ value.
This shows that the city block with the highest installed capacity is identical in both $q_{\mathrm{min}}$ scenarios and does not contain any well doublets with capacities between 1 and 5 [l/s].

\begin{table}[tb]
\scriptsize
\caption{Results of the six optimization scenarios for the city of Munich.}
\label{tab:results_01}
\centering
\begin{tabular}{llllllllllll}
\toprule
  \multicolumn{2}{c}{Scenario} & \multicolumn{3}{c}{Nr. of installed doublets} & \multicolumn{2}{c}{Doublet pumping rates [l/s]} & \multicolumn{2}{c}{Nr. of blocks} & \multicolumn{3}{c}{Block pumping rates [l/s]} \\
  \cmidrule(l){1-2}\cmidrule(lr){3-5}\cmidrule(lr){6-7}\cmidrule(lr){8-9}\cmidrule(lr){10-12}
 $q_{\mathrm{min}}$ [l/s] & $r_\Delta$ & Total & Max & Mean & Average Max & Mean & with & without & Total & Max & Mean \\
  &  &  &  &  & per block &  & doublets & doublets &  &  &  \\
\hline
1 & 1.5 & 24802 & 27 & 2.83 & 20.73 & 10.74 & 8232 & 519 & 274149.19 & 971.36 & 31.33 \\
   & 2 & 19134 & 20 & 2.19 & 21.28 & 12.05 & 8232 & 519 & 241056.50 & 762.76 & 27.55 \\
   & 3 & 14139 & 18 & 1.62 & 21.91 & 14.29 & 8232 & 519 & 212419.17 & 639.74 & 24.27 \\
\hline
5 & 1.5 & 14255 & 23 & 1.63 & 26.20 & 17.36 & 6342 & 2409 & 256201.37 & 971.36 & 29.28 \\ 
  & 2 & 11435 & 16 & 1.31 & 26.89 & 18.97 & 6342 & 2409 & 227090.83 & 762.76 & 25.95 \\
  & 3 & 8905 & 10 & 1.02 & 27.65 & 21.71 & 6342 & 2409 & 202045.02 & 639.74 & 23.09 \\
\bottomrule
\end{tabular}
\end{table}

\begin{figure}[tbh!]
\centering
\includegraphics[width=0.78\textwidth]{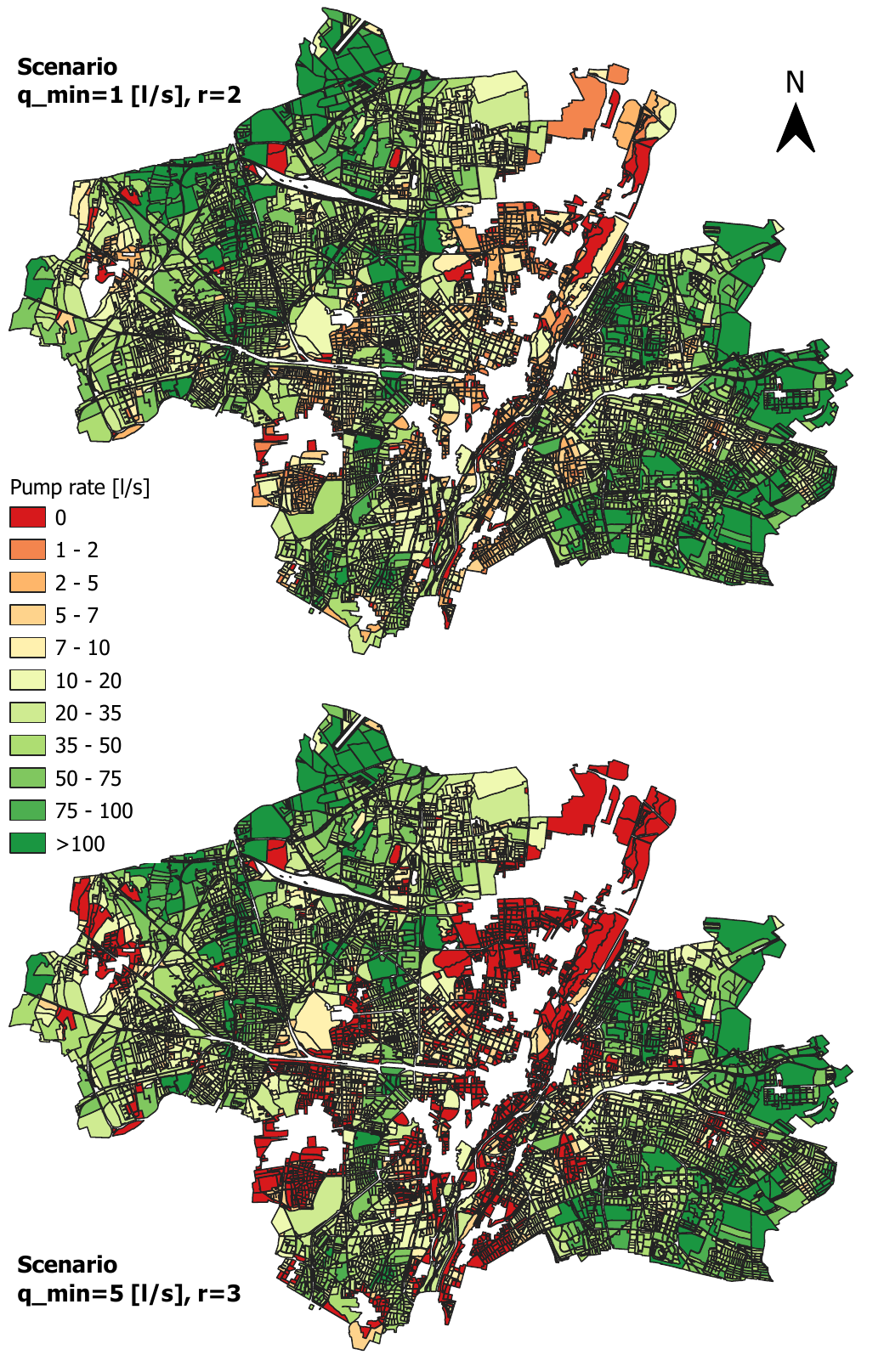}
\vspace{-.4cm}
\caption{Optimized pumping rates per city block for two scenarios $q_{\mathrm{min}}=1$ [l/s], $r_\Delta=2$ (top) and $q_{\mathrm{min}}=5$ [l/s], $r_\Delta=3$ (bottom).}
\label{fig:Munich_maps}
\end{figure}

The number of city blocks with and without installed doublets remains constant for scenarios with the same $q_{\mathrm{min}}$ value.
This is due to the fact that the parameter $r_\Delta$ controls the spacing between neighboring doublets and does not influence whether at least one single doublet is installed within a city block. 
From the scenarios with $q_{\mathrm{min}}=1$ [l/s] to $q_{\mathrm{min}}=5$ [l/s], the count of city blocks without installed doublets rises from 519 to 2409, respectively, of the total 8751 blocks considered. 
The difference of 1890 city blocks results from areas with low potential for thermal groundwater use mainly due to a lower groundwater thickness and corresponds to the blocks containing doublets with flow rates between 1 and 5 [l/s].
This observation is further evident in \autoref{fig:Munich_maps}, which shows the optimized potential of GWHP systems for two scenarios: $q_{\mathrm{min}}=1$ [l/s], $r_\Delta=2$ (top) and $q_{\mathrm{min}}=5$ [l/s], $r_\Delta=3$ (bottom).
The second scenario is more conservative, requiring more spacing between neighboring well doublets and considering only larger doublets, and this is also evident in the results (see Table~\ref{tab:results_01}).
Moreover, the results show that certain city regions exhibit an extensive potential for the thermal use of groundwater, making them especially well suited for the use of large multi-well GWHPs in 4GDH or 5GDHC systems. 
Conversely, areas in the inner city zone and around the city river Isar (see \autoref{fig:case_study}) with lower potential for thermal groundwater use (marked in red colors) are also observed to be less suitable for larger GWHPs.
These results are the consequence of unfavorable hydrogeologic conditions in the inner city zone, which were also observed in the original TAP publication \citep{Bottcher.2019}.

\begin{figure}[tb]
\centering
\includegraphics[width=0.85\textwidth]{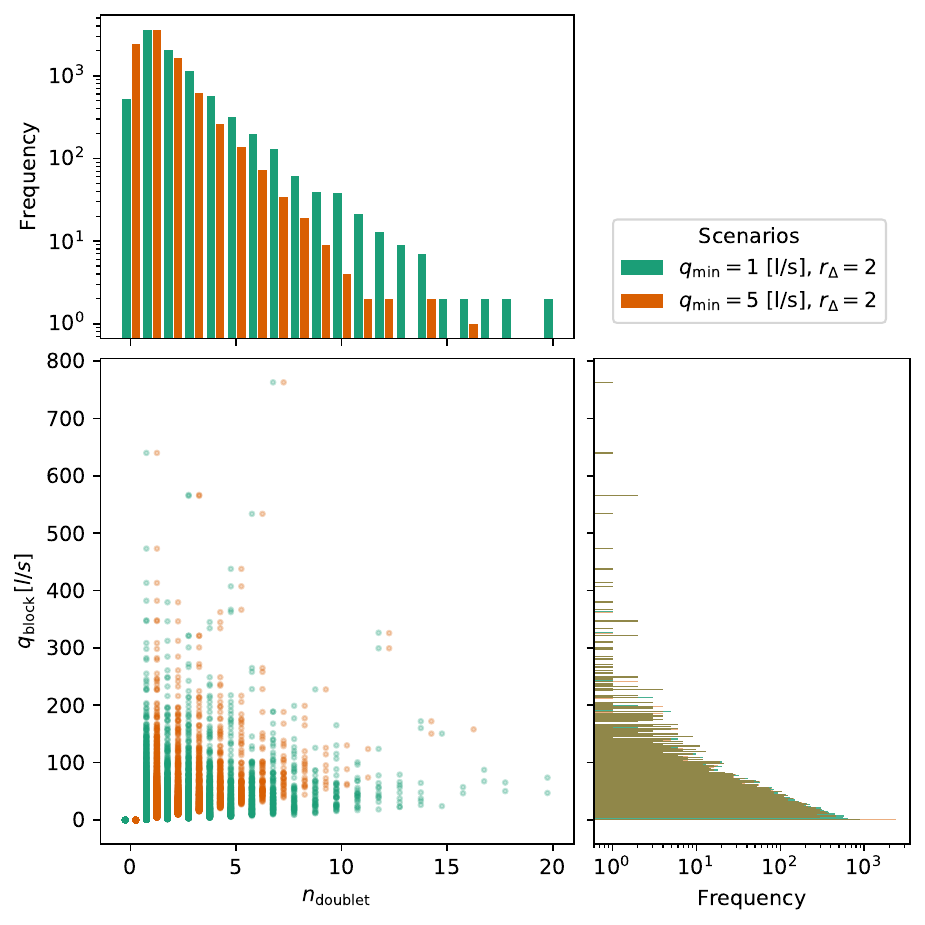}
\caption{Result comparison for two optimization scenarios.}
\label{fig:scenario_stats}
\end{figure}

\autoref{fig:scenario_stats} shows the statistical results for two optimization scenarios: $q_{\mathrm{min}}=1$ [l/s], $r_\Delta=2$ and $q_{\mathrm{min}}=5$ [l/s], $r_\Delta=2$.
The figure depicts the distributions of installed capacities $q_{\mathrm{block}}$ and the number of installed doublets $n_{\mathrm{doublet}}$ per city block, along with the correlation between these two variables. 
The distribution plots reveal that both parameters, $q_{\mathrm{block}}$ and $n_{\mathrm{doublet}}$, mostly follow an exponential distribution pattern.
This means that the frequency (count) of city blocks increases exponentially with decreasing capacity $q_{\mathrm{block}}$ and number of doublets $n_{\mathrm{doublet}}$ installed per block. 
Moreover, the number of blocks with only one installed doublet is the highest, followed by blocks with two or no doublets, depending on the scenario.
The scenario with $q_{\mathrm{min}}=5$ [l/s] contains more blocks without doublets compared to the scenario with $q_{\mathrm{min}}=1$ [l/s] because the first scenario excludes all blocks with only one doublet that has a capacity between 1 and 5 [l/s].
Similarly, due to the exclusion of smaller doublets in the first scenario, there are also fewer blocks in this scenario that contain a larger number of installed doublets (e.g. 10 doublets per block).
This is further supported by the central graph that illustrates the relationship between the number of doublets $n_{\mathrm{doublets}}$ in a block and its capacity $q_{\mathrm{block}}$.

\begin{figure}[tb]
\centering
\includegraphics[width=0.85\textwidth]{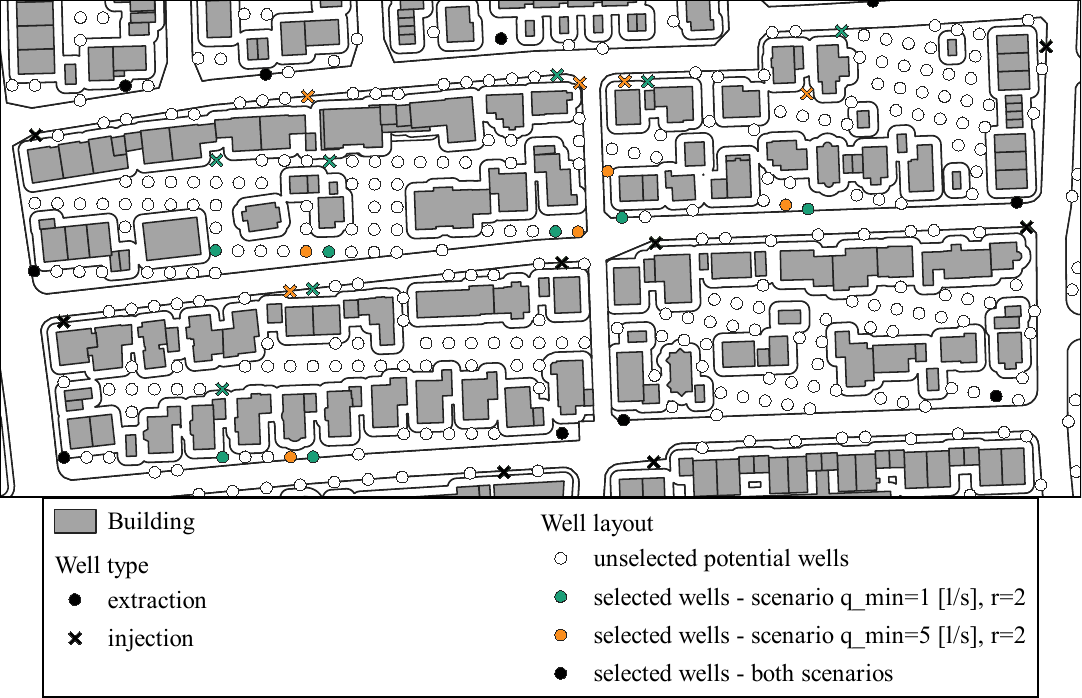}
\caption{Optimal well placement for two scenarios.}
\label{fig:result_wells}
\end{figure}

\autoref{fig:result_wells} presents an example of the optimal positioning of well doublets in four city blocks. 
The optimal well arrangements for two scenarios are depicted: $q_{\mathrm{min}}=1$ [l/s], $r_\Delta=2$ and $q_{\mathrm{min}}=5$ [l/s], $r_\Delta=2$. 
As discussed previously, only larger doublets are included in the second scenario. 
Consequently, certain city blocks in the second scenario have fewer but larger doublets compared to the first scenario, as can be seen in \autoref{fig:result_wells}.

\section{Discussion and Outlook}
\label{sec:discussion}

The presented optimization approach can effectively analyze the technical potential of thermal groundwater use by determining optimal arrangements of well doublets, their sizing, and well locations within the designated area. 
The approach can successfully quantify the maximum technically achievable groundwater pumping rate, and thus the exchanged thermal energy with the aquifer, taking into account relevant regulatory and hydraulic constraints.
Moreover, the considered optimization scenarios demonstrate the approach's versatility through the selection of specific optimization parameters. 
For example, by increasing the ratio $r_\Delta$, more conservative results are obtained, indicating reduced interaction between neighboring doublets. 
Similarly, selecting a larger value for the minimum capacity $q_{\mathrm{min}}$ focuses the analysis on large GWHP systems, which is particularly useful for investigating the potential for 4GDH and 5GDHC systems. 
The approach can also be used to study the potential of smaller systems, such as distributed GWHPs for individual households, by setting an upper threshold for the installed capacity of an individual well doublet. 
Therefore, the presented approach can serve as a valuable basis for both thermal groundwater management and urban energy planning. 
The new approach offers several advantages over existing methods, including simultaneous optimization of well placement and size, consideration of spatial heterogeneity of groundwater parameters, and hydraulic constraints for determining optimal pumping rates.
However, it has certain limitations, which are discussed in the following together with possible future improvements.

In this approach, the analytical formulas of the TAP method are incorporated into an optimization framework. 
As a result, certain limitations are inherited from the integration of the TAP method. 
The first one is that the wells of a doublet are fully aligned with the groundwater flow direction. 
In reality, this need not be the case, as wells can be placed in various configurations, such as in the corners of a city block, without strictly following the groundwater flow direction. 
To address this limitation, one of the first future improvements is to extend the approach to such cases. 
This extension involves using new analytical expressions applicable to well pairs that are not aligned with groundwater flow direction. 

The second limitation relates to having the same number of extraction and injection wells in a potential multi-well system. The proposed approach already provides flexibility in the design of multi-well systems because the doublets within a city block can be connected in different configurations. For example, multiple doublets can be combined into a larger system or several smaller single-doublet systems can be used. In practice, however, systems may have unbalanced combinations of extraction and injection wells due to hydrogeological conditions, such as one extraction well paired with two or more injection wells. Therefore, the approach can be further extended in the future to accommodate such scenarios.

It should be mentioned that the potential analysis study presented in this work did not consider existing systems within the city. 
Nevertheless, the approach presented is fully capable of including existing systems into the potential analysis. 
Including existing systems in the optimization problem \autoref{eq:Prob} is a straightforward process. 
One method is to set the optimization variables of the installed systems and wells to constant values using equality constraints. In particular, binary decision variables can be set to 1, representing the presence of the installed wells, while capacities can be set to their actual installed capacity values. Alternatively, optimization variables for existing systems need not be used, and their values can be substituted directly into the associated constraints with corresponding constant numerical values. 
By using either of these methods, the approach can efficiently account for existing systems and contribute to a more comprehensive analysis of thermal aquifer potential.

Furthermore, as stated in Section~\ref{sssec:well_selection}, the proposed approach does not directly consider thermal anomalies in the groundwater. 
However, certain aspects related to heat transport are addressed indirectly within the optimization approach. 
First, the prevention of internal hydraulic breakthrough within each system is achieved by applying the constraint \autoref{eq:cst_int_break}. 
Since thermal breakthrough normally occurs after hydraulic breakthrough, this constraint also ensures the prevention of internal thermal breakthrough \citep{Banks.2009}. 
Second, the constraints \autoref{eq:cst_ext_break} and \autoref{eq:cst_Dext_min} ensure sufficient spacing between neighboring well doublets, accounting for the hydraulic footprint of each doublet. 
Additionally, a potential doublet is defined on each line within a city block, and additional doublets are installed on parallel lines that are aligned with the groundwater flow direction. 
This geometric arrangement of potential doublets within a city block, along with the required spacing between neighboring doublets, results in a lower possibility of mutual thermal interference. 
Consequently, the new approach indirectly accounts for thermal effects such as thermal breakthroughs or negative thermal interactions between neighboring doublets. 
This only holds true for every city block individually without considering interactions with neighboring blocks.

On the other hand, when jointly optimizing multiple neighboring city blocks, it is crucial to consider the propagation of thermal plumes, since upstream systems can directly affect downstream ones. 
To account for this, an analytical model for estimating thermal plumes can be incorporated into the proposed optimization approach. 
One possible solution is to combine our approach with optimization concepts from the study by \citet{Halilovic.2023b}, where the spatial potential is optimized using the LAHM analytical model for thermal plume estimations. 
By incorporating these thermal aspects into the optimization process, the approach can be extended to the combined optimization of large areas, such as several neighboring city blocks or entire city districts. 
Moreover, there is the potential to combine the approach with energy system optimization models (ESOMs) used for optimal planning of urban energy systems.
For urban areas where thermal use of groundwater is a viable option, accurate representation of thermal potential in ESOMs is crucial \citep{Halilovic.2022}, especially in the context of optimal planning of future 4GDH and 5GDHC systems. 
In general, potential analysis methods based on analytical formulas are more suitable for integration into ESOMs because they are significantly less computationally demanding and complex compared to methods based on numerical groundwater simulations \citep{Halilovic.2022c}.

Finally, the approach can be extended to cost-related considerations, allowing for a holistic analysis that addresses both economic and environmental aspects of energy planning. 
This is particularly relevant for large groundwater uses with multiple wells for heating or cooling, as drilling costs become a significant factor. 
The challenge is to find the optimal balance between fewer, more expensive wells with larger diameters and multiple wells with smaller diameters. 
Additionally, the proposed approach can be fully integrated into GIS-based online management and energy planning tools, such as the web tool developed as part of the GEO.KW project \citep {Zosseder.2020b, EGU.2021}.

\section{Conclusion}\label{sec:Summary}

This paper presents a novel approach for determining the optimal sizing and placement of well doublets, with the overall goal of maximizing the technical potential of thermal groundwater use, i.e. the volume of pumped groundwater. 
The approach incorporates regulatory conditions, spatial variability of groundwater parameters, and key hydraulic constraints into the potential assessment process. 
The considered hydraulic constraints ensure acceptable drawdown levels in extraction wells and groundwater rise in injection wells, prevention of internal hydraulic breakthroughs, and adequate spacing between neighboring doublets. 
Analytic expressions describing these hydraulic constraints are integrated into a mixed-integer linear optimization problem allowing efficient application to relatively large areas.
The positioning of well doublets is based on the selection of predefined potential locations.

The proposed approach is applied to a real case study involving the optimization of the geothermal potential for each city block in Munich. 
To demonstrate the adaptability of the approach, six different optimization scenarios are used, differing in two parameters: the minimum capacity of a single installed doublet (1 and 5 [l/s]) and the external-internal well distance ratio (1.5, 2, and 3). 
The results prove the effectiveness and efficiency of the approach in identifying urban areas, or in this case city blocks, with favorable potential for large-scale GWHP systems as well as those unsuitable for such installations. 
Furthermore, the presented method provides comprehensive and optimized potential estimates that can be readily extended to other geographic locations. 
Thus, it is a valuable tool for thermal groundwater management and the integration of thermal groundwater potential into spatial energy planning, including the development of future 4GDH and 5GDHC systems. 
In addition, the method can provide valuable insights for well drilling and construction companies and housing associations. 
Finally, by coupling optimization techniques with potential analysis, the new method enables more thorough exploration and exploitation of the shallow geothermal potential, leading to an improved use of groundwater for heating and cooling purposes.

\section{Conflict of interest}
The authors declare no conflicts of interest. 

\section*{Acknowledgement}
\label{sec:Acknowledgement}
No external funding was received for conducting this study.

\bibliography{bibliography.bib}

\begin{thebibliography}{36}
\expandafter\ifx\csname natexlab\endcsname\relax\def\natexlab#1{#1}\fi
\providecommand{\url}[1]{\texttt{#1}}
\providecommand{\href}[2]{#2}
\providecommand{\path}[1]{#1}
\providecommand{\DOIprefix}{doi:}
\providecommand{\ArXivprefix}{arXiv:}
\providecommand{\URLprefix}{URL: }
\providecommand{\Pubmedprefix}{pmid:}
\providecommand{\doi}[1]{\href{http://dx.doi.org/#1}{\path{#1}}}
\providecommand{\Pubmed}[1]{\href{pmid:#1}{\path{#1}}}
\providecommand{\bibinfo}[2]{#2}
\ifx\xfnm\relax \def\xfnm[#1]{\unskip,\space#1}\fi
\bibitem[{Lund and Toth(2021)}]{Lund.2021}
\bibinfo{author}{J.~W. Lund}, \bibinfo{author}{A.~N. Toth},
\newblock \bibinfo{title}{Direct utilization of geothermal energy 2020 worldwide review},
\newblock \bibinfo{journal}{Geothermics} \bibinfo{volume}{90} (\bibinfo{year}{2021}) \bibinfo{pages}{101915}. \DOIprefix\doi{10.1016/j.geothermics.2020.101915}.
\bibitem[{Lund et~al.(2014)Lund, Werner, Wiltshire, Svendsen, Thorsen, Hvelplund, and Mathiesen}]{Lund.2014}
\bibinfo{author}{H.~Lund}, \bibinfo{author}{S.~Werner}, \bibinfo{author}{R.~Wiltshire}, \bibinfo{author}{S.~Svendsen}, \bibinfo{author}{J.~E. Thorsen}, \bibinfo{author}{F.~Hvelplund}, \bibinfo{author}{B.~V. Mathiesen},
\newblock \bibinfo{title}{{4th Generation District Heating (4GDH): Integrating smart thermal grids into future sustainable energy systems}},
\newblock \bibinfo{journal}{Energy} \bibinfo{volume}{68} (\bibinfo{year}{2014}) \bibinfo{pages}{1--11}. \DOIprefix\doi{10.1016/j.energy.2014.02.089}.
\bibitem[{Boesten et~al.(2019)Boesten, Ivens, Dekker, and Eijdems}]{Boesten.2019}
\bibinfo{author}{S.~Boesten}, \bibinfo{author}{W.~Ivens}, \bibinfo{author}{S.~C. Dekker}, \bibinfo{author}{H.~Eijdems},
\newblock \bibinfo{title}{5th generation district heating and cooling systems as a solution for renewable urban thermal energy supply},
\newblock \bibinfo{journal}{Advances in Geosciences} \bibinfo{volume}{49} (\bibinfo{year}{2019}) \bibinfo{pages}{129--136}. \DOIprefix\doi{10.5194/adgeo-49-129-2019}.
\bibitem[{Stauffer et~al.(2014)Stauffer, Bayer, Blum, Giraldo, and Kinzelbach}]{Stauffer.2014}
\bibinfo{author}{F.~Stauffer}, \bibinfo{author}{P.~Bayer}, \bibinfo{author}{P.~Blum}, \bibinfo{author}{N.~M. Giraldo}, \bibinfo{author}{W.~Kinzelbach}, \bibinfo{title}{Thermal use of shallow groundwater}, \bibinfo{publisher}{{CRC Press}}, \bibinfo{address}{Boca Raton, Fla.}, \bibinfo{year}{2014}.
\bibitem[{Banks(2012)}]{Banks.2012}
\bibinfo{author}{D.~Banks}, \bibinfo{title}{An introduction to thermogeology: ground source heating and cooling}, \bibinfo{publisher}{John Wiley \& Sons}, \bibinfo{year}{2012}.
\bibitem[{Lee et~al.(2006)Lee, Won, and Hahn}]{Lee.2006}
\bibinfo{author}{J.-Y. Lee}, \bibinfo{author}{J.-H. Won}, \bibinfo{author}{J.-S. Hahn},
\newblock \bibinfo{title}{Evaluation of hydrogeologic conditions for groundwater heat pumps: analysis with data from national groundwater monitoring stations},
\newblock \bibinfo{journal}{Geosciences Journal} \bibinfo{volume}{10} (\bibinfo{year}{2006}) \bibinfo{pages}{91--99}. \DOIprefix\doi{10.1007/BF02910336}.
\bibitem[{Busby et~al.(2009)Busby, Lewis, Reeves, and Lawley}]{Busby.2009}
\bibinfo{author}{J.~Busby}, \bibinfo{author}{M.~Lewis}, \bibinfo{author}{H.~Reeves}, \bibinfo{author}{R.~Lawley},
\newblock \bibinfo{title}{Initial geological considerations before installing ground source heat pump systems},
\newblock \bibinfo{journal}{Quarterly Journal of Engineering Geology and Hydrogeology} \bibinfo{volume}{42} (\bibinfo{year}{2009}) \bibinfo{pages}{295--306}. \DOIprefix\doi{10.1144/1470-9236/08-092}.
\bibitem[{Fry(2009)}]{Fry.2009}
\bibinfo{author}{V.~Fry},
\newblock \bibinfo{title}{{Lessons from London: regulation of open-loop ground source heat pumps in central London}},
\newblock \bibinfo{journal}{Quarterly Journal of Engineering Geology and Hydrogeology} \bibinfo{volume}{42} (\bibinfo{year}{2009}) \bibinfo{pages}{325--334}. \DOIprefix\doi{10.1144/1470-9236/08-087}.
\bibitem[{Schiel et~al.(2016)Schiel, Baume, Caruso, and Leopold}]{Schiel.2016}
\bibinfo{author}{K.~Schiel}, \bibinfo{author}{O.~Baume}, \bibinfo{author}{G.~Caruso}, \bibinfo{author}{U.~Leopold},
\newblock \bibinfo{title}{{GIS-based modelling of shallow geothermal energy potential for CO2 emission mitigation in urban areas}},
\newblock \bibinfo{journal}{Renewable Energy} \bibinfo{volume}{86} (\bibinfo{year}{2016}) \bibinfo{pages}{1023--1036}. \DOIprefix\doi{https://doi.org/10.1016/j.renene.2015.09.017}.
\bibitem[{Garc{\'\i}a~Gil et~al.(2022)Garc{\'\i}a~Gil, Garrido~Schneider, Mej{\'\i}as~Moreno, and Santamarta~Cerezal}]{GarciaGil.2022}
\bibinfo{author}{A.~Garc{\'\i}a~Gil}, \bibinfo{author}{E.~A. Garrido~Schneider}, \bibinfo{author}{M.~Mej{\'\i}as~Moreno}, \bibinfo{author}{J.~C. Santamarta~Cerezal},
\newblock \bibinfo{title}{Management and governance of shallow geothermal energy resources},
\newblock in: \bibinfo{booktitle}{Shallow Geothermal Energy}, \bibinfo{publisher}{Springer}, \bibinfo{year}{2022}, pp. \bibinfo{pages}{237--272}.
\bibitem[{Epting et~al.(2018)Epting, M{\"u}ller, Genske, and Huggenberger}]{Epting.2018}
\bibinfo{author}{J.~Epting}, \bibinfo{author}{M.~H. M{\"u}ller}, \bibinfo{author}{D.~Genske}, \bibinfo{author}{P.~Huggenberger},
\newblock \bibinfo{title}{Relating groundwater heat-potential to city-scale heat-demand: A theoretical consideration for urban groundwater resource management},
\newblock \bibinfo{journal}{Applied Energy} \bibinfo{volume}{228} (\bibinfo{year}{2018}) \bibinfo{pages}{1499--1505}. \DOIprefix\doi{10.1016/j.apenergy.2018.06.154}.
\bibitem[{Epting et~al.(2020)Epting, B{\"o}ttcher, Mueller, Garc{\'i}a-Gil, Zosseder, and Huggenberger}]{Epting.2020}
\bibinfo{author}{J.~Epting}, \bibinfo{author}{F.~B{\"o}ttcher}, \bibinfo{author}{M.~H. Mueller}, \bibinfo{author}{A.~Garc{\'i}a-Gil}, \bibinfo{author}{K.~Zosseder}, \bibinfo{author}{P.~Huggenberger},
\newblock \bibinfo{title}{City-scale solutions for the energy use of shallow urban subsurface resources -- bridging the gap between theoretical and technical potentials},
\newblock \bibinfo{journal}{Renewable Energy} \bibinfo{volume}{147} (\bibinfo{year}{2020}) \bibinfo{pages}{751--763}. \DOIprefix\doi{10.1016/j.renene.2019.09.021}.
\bibitem[{Casasso and Sethi(2017)}]{Casasso.2017}
\bibinfo{author}{A.~Casasso}, \bibinfo{author}{R.~Sethi},
\newblock \bibinfo{title}{{Assessment and mapping of the shallow geothermal potential in the province of Cuneo (Piedmont, NW Italy)}},
\newblock \bibinfo{journal}{Renewable Energy} \bibinfo{volume}{102} (\bibinfo{year}{2017}) \bibinfo{pages}{306--315}. \DOIprefix\doi{10.1016/j.renene.2016.10.045}.
\bibitem[{Muñoz et~al.(2015)Muñoz, Garat, Flores-Aqueveque, Vargas, Rebolledo, Sepúlveda, Daniele, Morata, and Ángel Parada}]{Munoz.2015}
\bibinfo{author}{M.~Muñoz}, \bibinfo{author}{P.~Garat}, \bibinfo{author}{V.~Flores-Aqueveque}, \bibinfo{author}{G.~Vargas}, \bibinfo{author}{S.~Rebolledo}, \bibinfo{author}{S.~Sepúlveda}, \bibinfo{author}{L.~Daniele}, \bibinfo{author}{D.~Morata}, \bibinfo{author}{M.~Ángel Parada},
\newblock \bibinfo{title}{{Estimating low-enthalpy geothermal energy potential for district heating in Santiago basin–Chile (33.5°S)}},
\newblock \bibinfo{journal}{Renewable Energy} \bibinfo{volume}{76} (\bibinfo{year}{2015}) \bibinfo{pages}{186--195}. \DOIprefix\doi{10.1016/j.renene.2014.11.019}.
\bibitem[{Garc{\'i}a-Gil et~al.(2015)Garc{\'i}a-Gil, V{\'a}zquez-Su{\~n}e, Alcaraz, Juan, S{\'a}nchez-Navarro, Montlle{\'o}, Rodr{\'i}guez, and Lao}]{GarciaGil.2015}
\bibinfo{author}{A.~Garc{\'i}a-Gil}, \bibinfo{author}{E.~V{\'a}zquez-Su{\~n}e}, \bibinfo{author}{M.~M. Alcaraz}, \bibinfo{author}{A.~S. Juan}, \bibinfo{author}{J.~{\'A}. S{\'a}nchez-Navarro}, \bibinfo{author}{M.~Montlle{\'o}}, \bibinfo{author}{G.~Rodr{\'i}guez}, \bibinfo{author}{J.~Lao},
\newblock \bibinfo{title}{Gis-supported mapping of low-temperature geothermal potential taking groundwater flow into account},
\newblock \bibinfo{journal}{Renewable Energy} \bibinfo{volume}{77} (\bibinfo{year}{2015}) \bibinfo{pages}{268--278}. \DOIprefix\doi{10.1016/j.renene.2014.11.096}.
\bibitem[{Pujol et~al.(2015)Pujol, Ricard, and Bolton}]{Pujol.2015}
\bibinfo{author}{M.~Pujol}, \bibinfo{author}{L.~P. Ricard}, \bibinfo{author}{G.~Bolton},
\newblock \bibinfo{title}{20 years of exploitation of the yarragadee aquifer in the perth basin of western australia for direct-use of geothermal heat},
\newblock \bibinfo{journal}{Geothermics} \bibinfo{volume}{57} (\bibinfo{year}{2015}) \bibinfo{pages}{39--55}. \DOIprefix\doi{10.1016/j.geothermics.2015.05.004}.
\bibitem[{Arola and Korkka-Niemi(2014)}]{Arola.2014}
\bibinfo{author}{T.~Arola}, \bibinfo{author}{K.~Korkka-Niemi},
\newblock \bibinfo{title}{The effect of urban heat islands on geothermal potential: examples from quaternary aquifers in finland},
\newblock \bibinfo{journal}{Hydrogeology Journal} \bibinfo{volume}{22} (\bibinfo{year}{2014}) \bibinfo{pages}{1953--1967}. \DOIprefix\doi{10.1007/s10040-014-1174-5}.
\bibitem[{Arola et~al.(2014)Arola, Eskola, Hellen, and Korkka-Niemi}]{Arola.2014b}
\bibinfo{author}{T.~Arola}, \bibinfo{author}{L.~Eskola}, \bibinfo{author}{J.~Hellen}, \bibinfo{author}{K.~Korkka-Niemi},
\newblock \bibinfo{title}{Mapping the low enthalpy geothermal potential of shallow quaternary aquifers in finland},
\newblock \bibinfo{journal}{Geothermal Energy} \bibinfo{volume}{2} (\bibinfo{year}{2014}) \bibinfo{pages}{1--20}. \DOIprefix\doi{10.1186/s40517-014-0009-x}.
\bibitem[{B{\"o}ttcher et~al.(2019)B{\"o}ttcher, Casasso, G{\"o}tzl, and Zosseder}]{Bottcher.2019}
\bibinfo{author}{F.~B{\"o}ttcher}, \bibinfo{author}{A.~Casasso}, \bibinfo{author}{G.~G{\"o}tzl}, \bibinfo{author}{K.~Zosseder},
\newblock \bibinfo{title}{Tap - thermal aquifer potential: A quantitative method to assess the spatial potential for the thermal use of groundwater},
\newblock \bibinfo{journal}{Renewable Energy} \bibinfo{volume}{142} (\bibinfo{year}{2019}) \bibinfo{pages}{85--95}. \DOIprefix\doi{10.1016/j.renene.2019.04.086}.
\bibitem[{Halilovic et~al.(2023)Halilovic, Böttcher, Zosseder, and Hamacher}]{Halilovic.2023a}
\bibinfo{author}{S.~Halilovic}, \bibinfo{author}{F.~Böttcher}, \bibinfo{author}{K.~Zosseder}, \bibinfo{author}{T.~Hamacher},
\newblock \bibinfo{title}{Optimization approaches for the design and operation of open-loop shallow geothermal systems},
\newblock \bibinfo{journal}{arXiv preprint arXiv:2307.11244}  (\bibinfo{year}{2023}). \DOIprefix\doi{10.48550/arXiv.2307.11244}.
\bibitem[{Park et~al.(2021)Park, Lee, Kaown, Lee, and Lee}]{Park.2021}
\bibinfo{author}{D.~Park}, \bibinfo{author}{E.~Lee}, \bibinfo{author}{D.~Kaown}, \bibinfo{author}{S.-S. Lee}, \bibinfo{author}{K.-K. Lee},
\newblock \bibinfo{title}{Determination of optimal well locations and pumping/injection rates for groundwater heat pump system},
\newblock \bibinfo{journal}{Geothermics} \bibinfo{volume}{92} (\bibinfo{year}{2021}) \bibinfo{pages}{102050}. \DOIprefix\doi{10.1016/j.geothermics.2021.102050}.
\bibitem[{Park et~al.(2020)Park, Kaown, and Lee}]{Park.2020}
\bibinfo{author}{D.~K. Park}, \bibinfo{author}{D.~Kaown}, \bibinfo{author}{K.-K. Lee},
\newblock \bibinfo{title}{Development of a simulation-optimization model for sustainable operation of groundwater heat pump system},
\newblock \bibinfo{journal}{Renewable Energy} \bibinfo{volume}{145} (\bibinfo{year}{2020}) \bibinfo{pages}{585--595}. \DOIprefix\doi{10.1016/j.renene.2019.06.039}.
\bibitem[{Halilovic et~al.(2022)Halilovic, B{\"o}ttcher, Kramer, Piggott, Zosseder, and Hamacher}]{Halilovic.2022a}
\bibinfo{author}{S.~Halilovic}, \bibinfo{author}{F.~B{\"o}ttcher}, \bibinfo{author}{S.~C. Kramer}, \bibinfo{author}{M.~D. Piggott}, \bibinfo{author}{K.~Zosseder}, \bibinfo{author}{T.~Hamacher},
\newblock \bibinfo{title}{Well layout optimization for groundwater heat pump systems using the adjoint approach},
\newblock \bibinfo{journal}{Energy Conversion and Management} \bibinfo{volume}{268} (\bibinfo{year}{2022}) \bibinfo{pages}{116033}. \DOIprefix\doi{10.1016/j.enconman.2022.116033}.
\bibitem[{Halilovic et~al.(2023)Halilovic, B{\"o}ttcher, Zosseder, and Hamacher}]{Halilovic.2023b}
\bibinfo{author}{S.~Halilovic}, \bibinfo{author}{F.~B{\"o}ttcher}, \bibinfo{author}{K.~Zosseder}, \bibinfo{author}{T.~Hamacher},
\newblock \bibinfo{title}{Optimizing the spatial arrangement of groundwater heat pumps and their well locations},
\newblock \bibinfo{journal}{Renewable Energy} \bibinfo{volume}{217} (\bibinfo{year}{2023}) \bibinfo{pages}{119148}. \DOIprefix\doi{10.1016/j.renene.2023.119148}.
\bibitem[{Kinzelbach(1987)}]{Kinzelbach.1987}
\bibinfo{author}{W.~Kinzelbach}, \bibinfo{title}{{Numerische Methoden zur Modellierung des Transports von Schadstoffen im Grundwasser}}, \bibinfo{address}{Oldenbourg}, \bibinfo{year}{1987}.
\bibitem[{Clyde and Madabhushi(1983)}]{Clyde.1983}
\bibinfo{author}{C.~G. Clyde}, \bibinfo{author}{G.~V. Madabhushi},
\newblock \bibinfo{title}{Spacing of wells for heat pumps},
\newblock \bibinfo{journal}{Journal of Water Resources Planning and Management} \bibinfo{volume}{109} (\bibinfo{year}{1983}) \bibinfo{pages}{203--212}. \DOIprefix\doi{10.1061/(ASCE)0733-9496(1983)109:3(203)}.
\bibitem[{Banks(2009)}]{Banks.2009}
\bibinfo{author}{D.~Banks},
\newblock \bibinfo{title}{Thermogeological assessment of open-loop well-doublet schemes: a review and synthesis of analytical approaches},
\newblock \bibinfo{journal}{Hydrogeology Journal} \bibinfo{volume}{17} (\bibinfo{year}{2009}) \bibinfo{pages}{1149--1155}. \DOIprefix\doi{10.1007/s10040-008-0427-6}.
\bibitem[{{Bayerisches Landesamt f{\"u}r Umwelt}(2012)}]{LfU.2012}
\bibinfo{author}{{Bayerisches Landesamt f{\"u}r Umwelt}}, \bibinfo{title}{{Planung und Erstellung von Erdw{\"a}rmesonden: LfU}}, \bibinfo{year}{2012}. \URLprefix \url{https://www.lfu.bayern.de/wasser/merkblattsammlung/teil3_grundwasser_und_boden/doc/nr_372.pdf}.
\bibitem[{Santos and Toffolo(2020)}]{Santos.2020}
\bibinfo{author}{H.~G. Santos}, \bibinfo{author}{T.~A. Toffolo}, \bibinfo{title}{Mixed integer linear programming with python}, \bibinfo{year}{2020}.
\bibitem[{Halilovic and Böttcher(2023)}]{Halilovic.2023c}
\bibinfo{author}{S.~Halilovic}, \bibinfo{author}{F.~Böttcher}, \bibinfo{title}{{Optimization of thermal groundwater use}}, \bibinfo{year}{2023}. \URLprefix \url{https://github.com/SHalilovic/Well-doublet-optimization}. \DOIprefix\doi{10.5281/zenodo.10031559}.
\bibitem[{Jordahl et~al.(2020)Jordahl, den Bossche, Fleischmann, Wasserman, McBride, Gerard, Tratner, Perry, Badaracco, Farmer, Hjelle, Snow, Cochran, Gillies, Culbertson, Bartos, Eubank, maxalbert, Bilogur, Rey, Ren, Arribas-Bel, Wasser, Wolf, Journois, Wilson, Greenhall, Holdgraf, Filipe, and Leblanc}]{kelsey_jordahl_2020_3946761}
\bibinfo{author}{K.~Jordahl}, \bibinfo{author}{J.~V. den Bossche}, \bibinfo{author}{M.~Fleischmann}, \bibinfo{author}{J.~Wasserman}, \bibinfo{author}{J.~McBride}, \bibinfo{author}{J.~Gerard}, \bibinfo{author}{J.~Tratner}, \bibinfo{author}{M.~Perry}, \bibinfo{author}{A.~G. Badaracco}, \bibinfo{author}{C.~Farmer}, \bibinfo{author}{G.~A. Hjelle}, \bibinfo{author}{A.~D. Snow}, \bibinfo{author}{M.~Cochran}, \bibinfo{author}{S.~Gillies}, \bibinfo{author}{L.~Culbertson}, \bibinfo{author}{M.~Bartos}, \bibinfo{author}{N.~Eubank}, \bibinfo{author}{maxalbert}, \bibinfo{author}{A.~Bilogur}, \bibinfo{author}{S.~Rey}, \bibinfo{author}{C.~Ren}, \bibinfo{author}{D.~Arribas-Bel}, \bibinfo{author}{L.~Wasser}, \bibinfo{author}{L.~J. Wolf}, \bibinfo{author}{M.~Journois}, \bibinfo{author}{J.~Wilson}, \bibinfo{author}{A.~Greenhall}, \bibinfo{author}{C.~Holdgraf}, \bibinfo{author}{Filipe}, \bibinfo{author}{F.~Leblanc}, \bibinfo{title}{geopandas/geopandas: v0.8.1}, \bibinfo{year}{2020}. \URLprefix
  \url{https://doi.org/10.5281/zenodo.3946761}. \DOIprefix\doi{10.5281/zenodo.3946761}.
\bibitem[{Zosseder et~al.(2022)Zosseder, Kerl, Albarr{\'a}n-Ord{\'a}s, Gossler, Kiecak, and Chavez-Kus}]{Zosseder.2022}
\bibinfo{author}{K.~Zosseder}, \bibinfo{author}{M.~Kerl}, \bibinfo{author}{A.~Albarr{\'a}n-Ord{\'a}s}, \bibinfo{author}{M.~Gossler}, \bibinfo{author}{A.~Kiecak}, \bibinfo{author}{L.~Chavez-Kus},
\newblock \bibinfo{title}{{Die hydraulischen Grundwasserverhältnisse des quartären und des oberflächennahen tertiären Grundwasserleiters im Großraum München}},
\newblock \bibinfo{journal}{{Geologica Bavarica}} \bibinfo{volume}{122} (\bibinfo{year}{2022}). \URLprefix \url{https://www.bestellen.bayern.de/shoplink/91122.htm}.
\bibitem[{Halilovic et~al.(2022{\natexlab{a}})Halilovic, Odersky, and Hamacher}]{Halilovic.2022}
\bibinfo{author}{S.~Halilovic}, \bibinfo{author}{L.~Odersky}, \bibinfo{author}{T.~Hamacher},
\newblock \bibinfo{title}{Integration of groundwater heat pumps into energy system optimization models},
\newblock \bibinfo{journal}{Energy} \bibinfo{volume}{238} (\bibinfo{year}{2022}{\natexlab{a}}) \bibinfo{pages}{121607}. \DOIprefix\doi{10.1016/j.energy.2021.121607}.
\bibitem[{Halilovic et~al.(2022{\natexlab{b}})Halilovic, Odersky, B{\"o}ttcher, Davis, Schulte, Zosseder, and Hamacher}]{Halilovic.2022c}
\bibinfo{author}{S.~Halilovic}, \bibinfo{author}{L.~Odersky}, \bibinfo{author}{F.~B{\"o}ttcher}, \bibinfo{author}{K.~Davis}, \bibinfo{author}{M.~Schulte}, \bibinfo{author}{K.~Zosseder}, \bibinfo{author}{T.~Hamacher},
\newblock \bibinfo{title}{Optimization of an energy system model coupled with a numerical hydrothermal groundwater simulation},
\newblock in: \bibinfo{booktitle}{Mapping the Energy Future-Voyage in Uncharted Territory-, 43rd IAEE International Conference, July 31-August 3, 2022}, \bibinfo{organization}{International Association for Energy Economics}, \bibinfo{year}{2022}{\natexlab{b}}.
\bibitem[{Zosseder et~al.(2022)Zosseder, Böttcher, Davis, Haas, Halilovic, Hamacher, Heller, Odersky, Pauw, Schramm, and and}]{Zosseder.2020b}
\bibinfo{author}{K.~Zosseder}, \bibinfo{author}{F.~Böttcher}, \bibinfo{author}{K.~Davis}, \bibinfo{author}{C.~Haas}, \bibinfo{author}{S.~Halilovic}, \bibinfo{author}{T.~Hamacher}, \bibinfo{author}{H.~Heller}, \bibinfo{author}{L.~Odersky}, \bibinfo{author}{V.~Pauw}, \bibinfo{author}{T.~Schramm}, \bibinfo{author}{S.~M. and}, \bibinfo{title}{{Schlussbericht zum Verbundprojekt GEO-KW: Kopplung des geothermischen Speicherpotenzials mit den wechselnden Anforderungen des urbanen Energiebedarfs zur effizienten Nutzung der regenerativen Energiequelle Grundwasser für die dezentrale Kälte- und Wärmebereitstellung in der Stadt}}, \bibinfo{type}{Technical Report}, Bundesministerium für Wirtschaft und Klimaschutz, \bibinfo{year}{2022}. \DOIprefix\doi{10.14459/2022md1692003}.
\bibitem[{B{\"o}ttcher et~al.(2021)B{\"o}ttcher, Davis, Halilovic, Odersky, Pauw, Schramm, and Zosseder}]{EGU.2021}
\bibinfo{author}{F.~B{\"o}ttcher}, \bibinfo{author}{K.~Davis}, \bibinfo{author}{S.~Halilovic}, \bibinfo{author}{L.~Odersky}, \bibinfo{author}{V.~Pauw}, \bibinfo{author}{T.~Schramm}, \bibinfo{author}{K.~Zosseder}, \bibinfo{title}{{Optimising the thermal use of groundwater for a decentralized heating and cooling supply in the city of Munich, Germany}}, \bibinfo{year}{2021}. \URLprefix \url{https://doi.org/10.5194/egusphere-egu21-14929}. \DOIprefix\doi{10.5194/egusphere-egu21-14929}.

\end{thebibliography}



\end{document}